\documentclass[12pt]{amsart}

\usepackage{fix-cm}
\usepackage[normalem]{ulem}

\usepackage{lineno}

\usepackage{soul}

\usepackage{graphicx}

\usepackage{amssymb, amsmath, amsthm, hyperref, euscript, color}

\usepackage[dvipsnames]{xcolor}

\usepackage{amscd}

\usepackage{tikz-cd}

\usepackage{bbm}

\usepackage{enumitem}

\usepackage[english]{babel}

\usepackage{mathtools}

\usepackage{makecell}

\usepackage[sorting=nty, giveninits=true, maxnames=10]{biblatex} 
\usepackage{csquotes}

\numberwithin{equation}{section}

\newtheoremstyle{theor}{6pt plus 1pt minus 1pt}{6pt plus 1pt minus 1pt}{\slshape}{}{\bfseries}{.}{5pt plus 1pt minus 1pt}{}
\newtheoremstyle{def}{6pt plus 1pt minus 1pt}{6pt plus 1pt minus 1pt}{}{}{\bfseries}{.}{5pt plus 1pt minus 1pt}{}
\newtheoremstyle{rmk}{6pt plus 1pt minus 1pt}{6pt plus 1pt minus 1pt}{}{}{\bfseries}{.}{5pt plus 1pt minus 1pt}{}
\newtheoremstyle{claim}{6pt plus 1pt minus 1pt}{6pt plus 1pt minus 1pt}{}{}{\bfseries}{.}{5pt plus 1pt minus 1pt}{}

\theoremstyle{theor}
\newtheorem{newstatement}{newstatement}
\newtheorem{lemma}[newstatement]{Lemma}
\newtheorem{theorem}[newstatement]{Theorem}
\newtheorem*{theorem*}{Theorem 2}

\newtheorem{proposition}[newstatement]{Proposition}

\theoremstyle{def}
\newtheorem{definition}[newstatement]{Definition}

\theoremstyle{rmk}
\newtheorem{remark}[newstatement]{Remark}

\newtheorem*{example*}{Example}

\theoremstyle{claim}

\theoremstyle{theor}
\newtheorem{thm}{Theorem}

\expandafter\let\expandafter\oldproof\csname\string\proof\endcsname
\let\oldendproof\endproof
\renewenvironment{proof}[1][\proofname]{%
  \oldproof[\slshape #1]%
}{\oldendproof}
\def\provedboxcontents#1{$\square$}

\makeatletter
\newsavebox\myboxA
\newsavebox\myboxB
\newlength\mylenA

\newcommand*\xoverline[2][0.75]{%
    \sbox{\myboxA}{$\m@th#2$}%
    \setbox\myboxB\null
    \ht\myboxB=\ht\myboxA%
    \dp\myboxB=\dp\myboxA%
    \wd\myboxB=#1\wd\myboxA
    \sbox\myboxB{$\m@th\overline{\copy\myboxB}$}
    \setlength\mylenA{\the\wd\myboxA}
    \addtolength\mylenA{-\the\wd\myboxB}%
    \ifdim\wd\myboxB<\wd\myboxA%
       \rlap{\hskip 0.5\mylenA\usebox\myboxB}{\usebox\myboxA}%
    \else
        \hskip -0.5\mylenA\rlap{\usebox\myboxA}{\hskip 0.5\mylenA\usebox\myboxB}%
    \fi}
\makeatother

\newcommand{\R}{\mathbb{R}}

\newcommand{\Z}{\mathbb{Z}}

\newcommand{\N}{\mathbb{N}}

\newcommand{\CP}{{\mathbb C\mkern-0.5mu\mathbb P}}

\newcommand{\RP}{\mathbb{RP}}

\newcommand{\simtimes }{\mathbin{\widetilde{\smash{\times}}}}

\DeclareMathOperator{\Pin}{Pin}

\makeatother

\begin{document}

\author[Sebasti\'an M. Camponovo]{Sebasti\'an M. Camponovo}

\title[Stably exotic non-orientable four-manifolds with small $\chi$]{Some stably exotic non-orientable four-manifolds with small Euler characteristic}

\address{Dipartimento di Matematica, Informatica e Geoscienze, Università di Trieste, via Weiss 8, 34128
Trieste (IT)}

\email{sebastianmatias.camponovo@phd.units.it}

\begin{abstract}
    We describe a cut-and-paste operation that produces new exotic manifolds starting from existing ones. Using the construction, we provide new examples of stably exotic manifolds. Our examples improve on previous work by achieving the smallest Euler characteristic currently known for stably exotic four-manifold pairs within their respective fundamental group classes. We produce examples with various infinite fundamental groups with torsion. Except for when the fundamental group is $\Z\times \Z/2$, the examples produced in this note minimize the Euler characteristic for non-orientable manifolds in their fundamental group class, providing new calculations for the non-orientable Hausmann-Weinberger invariant. We also calculate this invariant for some other non-cyclic groups, by deriving lower bounds from an observation on group (co)homology and exhibiting four-manifolds whose homeomorphism classes attain these bounds. 
\end{abstract}

\maketitle

\section{Introduction.}\label{Introduction} 
We say that two closed, smooth four-manifolds $M_1$ and $M_2$ are stably exotic if they are homeomorphic but not stably diffeomorphic, where stabilization refers to connected sum with copies of $S^2\times S^2$.
O\-ri\-en\-t\-a\-ble stably exotic manifolds do not exist, by a result of Gompf \cite{gompf1984stable}, on the other hand, work of Kreck \cite{kreck2006some} shows an abundance of examples in the non-orientable realm. The first contribution of this short note is to present a simple cut and paste
mechanism and to sample its efficiency by building new stably exotic four-manifolds with small Euler characteristic. The procedure can be summarized in
three steps.

\begin{enumerate}[label=\textbf{S.\arabic*}]
\item Begin with a pair of non-orientable homeomorphic four-manifolds $M_1$ and  $M_2$. Proceed by removing the tubular neighbourhood $\nu(\alpha_i)$ and $\nu(\alpha_i')$ of two disjoint locally flat orientation-reversing simple closed loops $\alpha_i,\alpha_i'\subset M_i$  to produce a compact four-manifold
\begin{equation*}
M_i^0\mathrel{\mathop:}= M_i\setminus (\nu(\alpha_i)\sqcup\nu(\alpha_i'))
\end{equation*}with boundary $\partial(M_i^0) = S^2\simtimes S^1\sqcup S^2\simtimes S^1$ for $i = 1, 2$. 

\item Assemble new four-manifolds\begin{equation}\label{Manifolds Produced}M_i(\alpha_i,\alpha_i';f_i)\mathrel{\mathop:}= M_i^0\cup_{f_i} ((S^2\simtimes S^1)\times [0,1])\end{equation} for $i = 1, 2$ where $f_i\colon S^2\simtimes S^1\sqcup S^2\simtimes S^1\to S^2\simtimes S^1\sqcup S^2\simtimes S^1$ is a homeomorphism.
Observe that (see Proposition \ref{homeocondition}) up to choosing homeomorphic manifolds $$h\colon M_1\to M_2$$ and as starting loops in our procedure $\alpha_i,\alpha_i'\subset M_i$ so that $h(\alpha_1)=\alpha_2$ and $h(\alpha_1')=\alpha_2'$ we can always choose $f_i$ for $i = 1, 2$ to obtain a homeomorphism $$M_1(\alpha_1,\alpha_1';f_1)\to M_2(\alpha_2,\alpha_2';f_2).$$
\item Find an invariant that distinguishes the diffeomorphism types
of the new four-manifolds $M_1(\alpha_1,\alpha_1';f_1)$ and $M_2(\alpha_2,\alpha_2';f_2)$.
\end{enumerate}

In many cases, the diffeomorphism invariant that distinguishes $M_1$ and $M_2$ is also capable of discerning the new manifolds $M_1(\alpha_1,\alpha_1';f_1)$ and $M_2(\alpha_2,\alpha_2';f_2)$. As a first product of this operation we obtain the following result.

\begin{thm}\label{Theorem A}
    Let $p\in\Z_{>0}$ be odd. There are closed smooth non-orientable four-manifolds, $X_{2p}$ and $Y_{2p}$ that satisfy the following properties: 
    \begin{enumerate}
        \item There is a homeomorphism between $X_{2p}$ and $Y_{2p}.$
        \item $\chi(X_{2p})=1,\;w_2(X_{2p})=0$ and $\pi_1(X_{2p})=\Z\times\Z/2p$. 
        \item  There is no diffeomorphism between $X_{2p}\#(n-1)(S^2\times S^2)$ and $Y_{2p}\#(n-1)(S^2\times S^2)$ for any $n\in\Z_{>0}.$
        \item There is a diffeomorphism $$\widehat{X}_{2p}\to\widehat{Y}_{2p}$$ between the orientation 2-covers.
        \item There is a diffeomorphism $$X_{2p}\#\CP^2\to Y_{2p}\#\CP^2.$$
    \end{enumerate}
\end{thm}

Using a similar cut-and-paste procedure devised by Bais-Torres \cite{bais} we produce others families of stably exotic manifolds.

\begin{thm}\label{Theorem A2}
Let $p,q_1,\dots,q_m\in\Z_{>1}$ be odd. There are closed smooth non-orientable four-manifolds, $A=A(2p)(q_1,\dots,q_m)$ and $B=B(2p)(q_1,\dots,q_m)$ that satisfy the following properties: 
    \begin{enumerate}
        \item There is a homeomorphism between $A$ and $B$ 
        \item $\chi(A)=1,\;w_2(A)=0$ and \[\pi_1(A)=\Z/2\times(\Z/p*\Z/q_1*\dots*\Z/q_m)\]. 
        \item  There is no diffeomorphism between $A\#(n-1)(S^2\times S^2)$ and $B\#(n-1)(S^2\times S^2)$ for any $n\in\Z_{>0}.$
        \item There is a diffeomorphism $$\widehat{A}\to\widehat{B}$$ between the orientation 2-covers.
        \item There is a diffeomorphism $$A\#\CP^2\to B\#\CP^2.$$
   \end{enumerate}
\end{thm}

Examples of stably exotic four-manifolds with minimizing Euler characteristic in their fundamental group class were produced by Cappell-Shaneson \cite{cappell1976some} and Fintushel-Stern \cite{fintushel1981exotic} for $G=\Z/2$, Bais-Torres \cite{bais} for  $G=\Z/2p$ and Torres \cite{torres2017smooth} for $G=\Z/2*\Z/2$. Moreover, none of these examples need any connected sum with additional copies of $S^2\times S^2$ to pin down the homeomorphism class of the manifold realizing the required fundamental group. For $p>1$ the manifolds of Theorem \ref{Theorem A} and Theorem \ref{Theorem A2} belong to this class of examples, (see Proposition \ref{minimize} and Proposition \ref{secondminimality}). 
Another application of our procedure is to produce non-orientable manifolds that minimize the Euler characteristic in their fundamental group class. Recall that Hausmann-Weinberger \cite{Hausmann} defined an invariant of a finitely presented group $G$ as the minimizer of the Euler characteristic of closed oriented four-manifolds with fundamental group $G$. In general, the invariant is hard to compute, and unknown except for few cases, see \cite{Hausmann,kirkhausmann, kotschick}. Recently, Torres \cite[Definition D]{torres2026topologicalclassificationcertainnonorientable} proposed the analogous definition of an invariant $q(G)$ in the non-orientable settings and computed it for some finite cyclic groups. The content of our second theorem is to enlarge the family of groups $G$ for which $q(G)$ is known, through the cut-and-paste techniques introduced in this note and in \cite{bais}.  

\begin{thm}\label{Theorem B} Let $p,q_1,\dots,q_m\in\Z_{>1}$ be odd:
\begin{enumerate}
    \item  $q(\Z\times \Z/2p)=1$.
    \item $q(\Z/2\times(\Z/p*\Z/q_1*\dots*\Z/q_m))=1$.
    \item $q(\Z/2\times(\Z/2*\Z/q_1*\dots*\Z/q_m))=1$.
\end{enumerate}
\end{thm}

A simple use of covering space theory allows us to prove the following. 

\begin{proposition}\label{propcoveringq(G)}
    Let $G$ be an infinite virtually cyclic group. Then every closed non-orientable four-manifold with $\pi_1(M)=G$ has $\chi(M)\geq0.$ In particular $q(G)\geq0$.
\end{proposition}

From here we can compute the invariant $q(G)$ for infinitely many examples of groups that can be realized as fundamental group of four-dimensional mapping torus.

\begin{thm}\label{Theorem C}
     Let $M$ be a closed orientable three-manifold with finite fundamental group that admits an orientation-reversing self-diffeomorphism $f\colon M\to M$. Consider the mapping torus \[M_f=\frac{M\times[0,1]}{(x,0)\sim(f(x),1)}.\]
        Denote by $\Gamma_f=\pi_1(M_f)$, then $q(\Gamma_f)=0$.
\end{thm}

\begin{remark}\label{nonorientablelens}
    As an application of Theorem \ref{Theorem C} we immediately have  \[q(\Z/2\times \Z)=0\] realized by $M=\RP^3\simtimes S^1.$ 
    
    It is a classic result from Bonahon \cite[Theoreme 3]{Bonahon} that the lens space $L(p,q)$ admits an orientation-reversing self-diffeomorphism $\overline{\varphi}\colon L(p,q)\to L(p,q)$, if and only if the pair $(p,q)$ satisfies the following arithmetic condition \begin{equation}\label{lporientation}
    q^2\equiv-1\pmod p.
\end{equation}
When equation \eqref{lporientation} is satisfied, we define the mapping torus \begin{equation}\label{mappingtorus}
    L(p,q)\simtimes S^1=\frac{L(p,q)\times [0,1]}{(x,0)\sim(\overline{\varphi}(x),1)}.
\end{equation} The mapping torus structure and the work of Bonahon give us an explicit presentation of the fundamental group 
\[\Gamma_q=\pi_1(L(p,q)\simtimes S^1)=
\left\langle x,y \,\middle|\,
x^p,\;yxy^{-1}x^{-q}
\right\rangle,\; q^2\equiv -1\pmod p,\] therefore \[q(\Gamma_q)=0.\] 
\end{remark}

We also show that, in general, $q(G)$ can take arbitrarily negative values. 

\begin{thm}\label{Theorem D}
    Let $\Sigma_g$ be a closed orientable surface of genus $g\geq1$. Let $\Gamma_g=\pi_1(\Sigma_g)\times\Z/2$, then \[q(\Gamma_g)=2-2g,\] realized by $M_g=\Sigma_g\times\RP^2$.\\
    Let $p,q,r\in\Z_{>0}$ with $p$ even, $q^2\equiv-1\pmod p$ and $r$ odd. Let $\Gamma_{r}$ be the group \[\Gamma_{r}=\langle \alpha,\beta,\gamma|\;\beta^p,\alpha\beta\alpha^{-1}\beta^{-q}, y^{2},\alpha\gamma^{2r}\rangle=\Z/2\times\Z/2r,\] then \[q(\Gamma_{r})=1,\]
    realized by $M_{p,q,r}=(L(p,q)\simtimes S^1)\#_{S^1}\RP^4$.
\end{thm}
\begin{remark}
    Not every even integer $p>1$ admits an integer $q$ satisfying the condition in Theorem \ref{Theorem D}. For instance, $(p,q)=(2,1)$ is an admissible pair, whereas no such $q$ exists for $p=6$.
\end{remark}
We note that, recently, Torres \cite[Theorem B]{torres2026novelconstructionhomotopy4spheres} produced a myriad of examples of exotic manifolds in the homeomorphism class of $\RP^4$. Although it is not yet known whether the set contains more than two elements that are pairwise non-diffeomorphic, all of the manifolds produced by Torres can be used as inputs for the construction presented in this note, yielding corresponding relative versions of Theorem \ref{Theorem A} and \ref{Theorem A2}. 

The note is organized as follows. Section \ref{sectionpin} contains a brief description of Pin$^+$-structure on four-manifolds and the $\eta$-invariant of the twisted
Dirac operator, this will be the key invariant that we use to show that the manifolds appearing in Theorem \ref{Theorem A} are stably exotic. In Section \ref{cutandpasteprocedure} we present the construction of self circle-sum along loops in a non-orientable four-manifold and establish some interesting features of the mechanism. In Section \ref{sec:circlesum} we recall the procedure of circle-sum(s) already present in the literature and lay down the strategy blueprint that we will follow to construct the examples of Theorem \ref{Theorem A2} and Theorem \ref{Theorem B}. In Section \ref{bb} we present the building blocks manifolds that will be the inputs of our procedures. In Section \ref{constructA} we finally produce examples of manifolds that satisfy the conditions stated in Theorem \ref{Theorem A}.
We will verify the latter in Section \ref{certificateA}. Section \ref{sec:a2} and Section \ref{sec:proofA2} are devoted to the assemblage and study of the manifolds in Theorem \ref{Theorem A2}.
After this we move on to some calculations of the non-orientable Hausmann-Weinberger in Section \ref{HWinvariant}. Finally, in Appendix \ref{deficiency}, through the use of group homology, we compute the rank of the groups claimed in the various proposition. \subsection{Acknowledgements}We thank Rafael Torres for suggesting this project and for the constant support and guidance in its development. We are grateful to Marc Kegel for identifying an error in an earlier draft of this note. S.M.C. is a member of GNSAGA, Istituto
Nazionale di Alta Mate\-matica ‘Francesco Severi’, Italy.

\subsection*{Tools used} AI-assisted tools were employed in the final stage of manuscript preparation to detect grammatical errors, notational inconsistencies, and minor mathematical errors

\section{Pin$^+$-structures and a spectral invariant to distinguish smooth structures}\label{sectionpin} In this section, we establish the
conventions and notations on Pin$^+$-structure on a four-manifold that are
used in the sequel. The reader is directed towards \cite{kirby1990pin, stolz1988exotic} for further details. Since it will be used to distinguish the diffeomorphism type of our
non-orientable four-manifolds we recall some fundamental results of the theory of the $\eta$-invariant of a closed smooth four-manifold that is equipped with a Pin$^+$-structure, see \cite{stolz1988exotic} for background on this invariant. 

A Pin$^+$-structure $\varphi_M$ on a smooth four-manifold
$M$ will be denoted by $(M, \varphi_M )$. In the sequel we will also refer to the pair $(M, \varphi_M)$ as
a Pin$^+$-manifold. For any such four-manifold $M$ with first cohomology
group $H^1(M; \Z/2) =\Z/2$, there are two Pin$^+$-structures $\{(M, \pm\varphi_M )\}$,
and their classes $\{[(M, \pm\varphi_M )]\}$ are mutual inverses in the fourth Pin$^+$-
bordism group $\Omega^{\textnormal{Pin}^+}_4=\Z/16$ \cite[Corollary 6.4]{stolz1988exotic}, \cite[p. 184]{kirby1990pin}. In the general case the cohomology group $H^1(M;\Z/2)$ acts transitively and effectively on the set of Pin$^+$-structures, parameterizing such structures \cite[Corollary 6.4]{stolz1988exotic}, \cite[p. 184]{kirby1990pin}. We remark that we will distinguish the diffeomorphism classes by looking at the Pin$^+$-structures
on the tangent bundle of our four-manifolds through the $\eta$-invariant.

We now state results related to Pin$^+$-structure and the $\eta$-invariant that will be paramount for our purposes.

\begin{theorem}\label{Theorem NonOrientable 3}\

$\bullet$  \textup{(Kirby-Taylor) \cite[Theorem 5.2]{kirby1990pin}}. There is a group isomorphism \[\Omega^{\Pin^+}_4\rightarrow \Z/16\] and the class of the real projective four-space\begin{equation*}[(\RP^4, \pm \varphi_{\RP^4})] = \pm 1\in \Omega^{\Pin^+}_4\cong \Z/16.\end{equation*}is a generator.

$\bullet$  \textup{(Stolz) \cite[Proposition 4.3]{stolz1988exotic}}. The $\eta$-invariant$\pmod {2\Z}$ is a complete $\Pin^+$-cobordism invariant and it is additive under disjoint unions. 
\end{theorem}

The behavior of exoticity and stability in the family of manifolds that we will construct will be controlled by means of the following theorem.

\begin{theorem}\label{torresexotic}\textup{(Torres) \cite[Theorem 2]{torres2025exotic}}. Let $M_1$ and $M_2$ be closed non-orientable smooth four-manifolds that admit a \textnormal{Pin}$^+$-structure and that have the same Euler
characteristic. \\If for every \textnormal{Pin}$^+$-structure ${(M_i, \varphi_{M_i} ) : i = 1, 2}$, the
values of the corresponding $\eta$-invariants satisfy
\begin{equation}
    \eta(M_1, \varphi_{M_1})\neq \eta(M_2, \varphi_{M_2} )\; \textnormal{mod}\;2\Z,
\end{equation}
then there is no diffeomorphism
\begin{equation}
    M_1\# n(S^2\times S^2)\to M_2 \# n(S^2\times S^2)
\end{equation}

for any $n \in\N$.
\end{theorem}
\section{Building some non-orientable four-manifolds.}\label{lego}
\subsection{Cut-and-paste procedures.}\label{cutandpasteprocedure}We introduce a topological cut-and-paste procedure to construct a non-orientable manifold $X$ from a given initial one $Y$.
Our goal is to study relations between the topological invariants of $X$ and $Y$ and show that, in suitable cases, if one performs the operation starting from $Y_1$ and $Y_2$ homeomorphic but not diffeomorphic manifolds the resulting $X_1$ and $X_2$ will share the same property.

We recall that given an orientation-reversing locally flat simple loop $\alpha$ in a topological four-manifold its tubular neighbourhood is homeomorphic to the total space of the non-orientable three-disk bundle over the circle $\nu(\alpha)=D^3\simtimes S^1$ with boundary the non-orientable two-sphere bundle over $S^1$; $\partial\nu(\alpha)=S^2\simtimes S^1$. For details on the existence and uniqueness of tubular neighbourhood in the topological category the reader is directed to \cite[Section 9]{freedman2014topology},\cite[Chapter 5]{friedl2024foundations}. 

We proceed now to give a description of the procedure. Given a non-orientable four-manifold $M$ consider disjoint locally flat simple loops $\alpha,\alpha'\subset M$ such that the orientation character $w(\alpha)=w(\alpha')\neq0$ and denote by \begin{equation}\label{svk}
    M(\alpha,\alpha';f):=M \setminus(\nu(\alpha) \sqcup\nu(\alpha'))\cup_f((S^2\simtimes S^1)\times [0,1])
\end{equation}
where $f\colon S^2\simtimes S^1\sqcup S^2\simtimes S^1\to S^2\simtimes S^1\sqcup S^2\simtimes S^1$ is a homeomorphism.

The gluing map $f$ decomposes naturally into two maps defined on each connected component of the boundary; $f=(f^1,f^2)$ where $f^i$ is a map in the mapping class group of $S^2\simtimes S^1$. Here, and in the following, $f$ is always assumed to preserve the order of the boundary components.
In order to pin down the topological features of $M(\alpha,\alpha';f)$, we need to understand how the isotopy class of $f$ influences our procedure.
\begin{theorem}\label{isotopyclass}
\textup{(Kim-Raymond)\cite[Theorem 2.3]{kim1990diffeotopy} (Miller-Naylor)\cite[Remark 3.12]{miller}}. The diffeotopy group of $S^2\simtimes S^1$ is isomorphic to $\Z/2\times\Z/2$. It is generated by a reflection $\rho$, and a twist about a 2-sphere fiber, $\tau$.
\end{theorem}
We remark that our construction alters the fundamental group, but preserves many topological data of our original manifolds, such as the Euler characteristic $\chi(M)$ and the Kirby-Siebenmann invariant $\text{KS}(M)$. 

\begin{lemma}\label{costruzione}
    Let $M$ be a closed, non-orientable, four-manifold. Fix a presentation $\langle G\;|\;R\rangle $ of $\pi_1(M)$.\\Denote by $M(\alpha,\alpha';f)$ the manifold obtained following the construction above. Then the following holds: 
    \begin{enumerate}
        \item $\chi(M)=\chi(M(\alpha,\alpha';f)).$
        \item $\pi_1(M(\alpha,\alpha';f))=\langle G,s,t\;|\;R,\alpha^{\epsilon}s, tst^{-1}(\alpha')^{\epsilon'}\rangle$ with $\epsilon,\epsilon'\in\{\pm1\}$.
        \item There exists a 5-manifold $W$ such that $V=(W;M,M(\alpha,\alpha';f))$ is a cobordism. 
        \item \textnormal{KS}(M)=\textnormal{KS}$(M(\alpha,\alpha';f))$.
        \item Moreover, if $M$ is a smooth four-manifold and admits a \textnormal{Pin}$^+$\textnormal{-}structure, we can always find an $f$ such that $W$ is a \textnormal{Pin}$^+$\textnormal{-}cobordism. In this case we have $w_2(M)=0$ if and only if $w_2(M(\alpha,\alpha';f))=0$.
    \end{enumerate}
\end{lemma}
\begin{proof}
    The first clause of the lemma follows from (co)homological computations through Mayer-Vietoris as we have: $$\chi(M(\alpha,\alpha';f))=\chi(M)+\chi(S^2\simtimes S^1\times I)-2\chi(S^2\simtimes S^1)=\chi(M).$$ 
    For the second Item we must carefully track the data of our identification. We fix a presentation of the fundamental group of $M$ $$\pi_1(M)=\langle G|R\rangle$$ where $G$ is a finite set of generators and $R$ is a finite set of relators. We first observe that because of codimension and a transversality argument we have $\pi_1(M\setminus(\nu(\alpha) \sqcup\nu(\alpha')))=\pi_1(M),$ while it follows from the long exact sequence of a fibration that $\pi_1(S^2\simtimes S^1)\cong\Z.$ Since $[0,1]$ is contractible we have that $S^2\simtimes S^1\times [0,1]$ is homotopy equivalent to $S^2\simtimes S^1$ and we fix a generator $\langle s\rangle=\pi_1(S^2\simtimes S^1\times[0,1])$. Let $t\subset M(\alpha,\alpha';f)$ be the loop obtained by joining the ends of the $[0,1]$ factor through an arc joining $\partial\nu(\alpha)$ and $\partial\nu(\alpha')$. Apply the Seifert-van Kampen theorem to \ref{svk} to obtain the aforementioned presentation of the fundamental group, where the exponents $\epsilon,\epsilon'$ in the relators set depend only on the choice of $f=(f^1,f^2)$. Notice that all four possibilities are realizable, via Theorem \ref{isotopyclass}. For the cobordism in the third clause we start by considering the trivial cylinder $M\times[0,1]$. Remove the copy of the loops in the upper lid $M\times\{1\}$ and identify with the boundary in $S^2\simtimes S^1\times[0,1]$ via $f$ \begin{equation}\label{cobordismglue}(M\times\{1\})\setminus((\nu(\alpha)\sqcup \nu(\alpha'))\times\{1\})\cup_f(S^2\simtimes S^1\times[0,1])
    \end{equation} obtaining the desired cobordism between $M$ and $M(\alpha,\alpha';f)$ (see \ref{svk}). Notice that the operation described in \ref{cobordismglue} corresponds to the gluing of a 5-dimensional round 1-handle \cite[Definition 7.8]{akbulut20164} to the top boundary of the trivial cylinder $M\times[0,1]$. 
    Hsu \cite{Hsu} and Lashof-Taylor \cite{Lashof} have shown that the Kirby-Siebenmann invariant is a cobordism
invariant, this and the third Item of the Lemma imply $\textnormal{KS}(M)=\textnormal{KS}(M(\alpha,\alpha';f))$ for every choice of $f$.

    When $(M,\varphi_M)$ is already a smooth Pin$^+$-manifold we want to find a Pin$^+$-structure on $(W,\varphi_W)$ that extends the one on $M$. Notice that a Pin$^+$-structure on $(M,\varphi_M)$ induces a Pin$^+$-structure on $M\setminus(\nu(\alpha) \sqcup\nu(\alpha'))$ as a codimension zero submanifold \cite[Construction 1.5]{kirby1990pin}, that restricts to Pin$^+$-structures on each connected component of $\partial(M\setminus(\nu(\alpha) \sqcup\nu(\alpha')))=S^2\simtimes S^1\sqcup S^2\simtimes S^1$. 
    
    To ease notation, we denote by $\partial_i=S^2\simtimes S^1,$ for $i=1,2$, the boundary component coming respectively from $\alpha,\alpha'$.

    \begin{equation}
        (\partial_i(M\setminus(\nu(\alpha) \sqcup\nu(\alpha'))),\varphi_M^i)=(S^2\simtimes S^1,\varphi_{\partial_i}),
    \end{equation}
Notice that given a Pin$^+$-structure on $(M,\varphi_M)$ we can trivially extend it to a Pin$^+$-structure on the cylinder $(M\times I,\varphi_{M\times I})$. This again induces a Pin$^+$-structure on the upper-lid boundary of our cobordism, and consequently induces a Pin$^+$-structure on $M\times\{1\}\setminus(\nu(\alpha)\times{\{1\}} \sqcup\nu(\alpha')\times{\{1\}})$ when we remove the copy of the loops $\alpha$ and $\alpha'$ in the slice $M\times\{1\}$. When we glue the remaining piece of the cobordism we can always match the Pin$^+$-structure  $\varphi_{\partial_i}$ by post-composing \cite{Aitchison}, if needed, one of the components of $f=(f^1,f^2)$ by $\tau$. The claim $w_2(M)=0$ if and only if $w_2(M(\alpha,\alpha';f))=0$ follows from the discussion above concerning the existence of the Pin$^+$-cobordism, together with the characterization of the existence of Pin$^+$-structures in terms of the vanishing of $w_2$ \cite[Lemma 1.3]{kirby1990pin}.
\end{proof}

Given its paramount role for our goal, we single out the following homeomorphism criterion of our procedure. 

\begin{proposition}\label{homeocondition}
    Let $M_1$ and $M_2$ be topological four-manifolds and let $$\{\alpha_i,\alpha_i'\subset M_i:i=1,2\}$$ be disjoint locally flat orientation-reversing simple loops. Suppose that there is a homeomorphism $$h\colon M_1\setminus(\nu(\alpha_1)\sqcup\nu(\alpha_1'))\to M_2\setminus(\nu(\alpha_2)\sqcup\nu(\alpha_2')).$$ For every choice of $f_1$ as in (\ref{svk}) we can always find a map $f_2$ such that the four-manifolds $M_1(\alpha_1,\alpha_1';f_1)$ and $M_2(\alpha_2,\alpha_2';f_2)$ are homeomorphic. 
\end{proposition}
 \begin{proof} Choose a presentation of the fundamental group and a gluing map $$(f_1^1,f_1^2)=f_1\colon \partial(M_1\setminus(\nu(\alpha_1)\sqcup\nu(\alpha_1')))\to S^2\simtimes S^1\sqcup S^2\simtimes S^1,$$ as described in Lemma \ref{costruzione}. The homeomorphism $$h:M_1\setminus(\nu(\alpha_1)\sqcup\nu(\alpha_1'))\to M_2\setminus(\nu(\alpha_2)\sqcup\nu(\alpha_2'))$$ restricts to a homeomorphism of the boundary, without loss of generality, up to reordering the subscript, suppose that $h$ maps $$\partial\nu(\alpha_1)\to \partial\nu(\alpha_2) \text{ and } \partial\nu(\alpha_1')\to \partial\nu(\alpha_2').$$ We proceed by defining $f_2=(f_2^1,f_2^2)$ as $$f_2^1=f_1^1\circ h^{-1} \text{ and }f_2^2= f_1^2\circ h^{-1},$$ and it follows from the decomposition \begin{equation}
    M_i(\alpha_i,\alpha_i';f_i):=M_i \setminus(\nu(\alpha_i) \sqcup\nu(\alpha_i'))\cup_{f_i}((S^2\simtimes S^1)\times [0,1]),
\end{equation}
 that the identity map $$id\colon((S^2\simtimes S^1)\times [0,1])\to ((S^2\simtimes S^1)\times [0,1]),$$ extends the original homeomorphism $h$. \end{proof}
\begin{remark}\label{pinremark} 
    An analogue of Proposition \ref{homeocondition} carries over with the same proof to the smooth category. More care is needed when we try to achieve extensions and compare Pin$^+$-structures with those on the original Pin$^+$-manifold $(M,\varphi)$ as described in Lemma \ref{costruzione}. Recall that for a smooth Pin$^+$-manifold $M$ that has $H^1(M;\Z/2)\cong \Z/2$ its Pin$^+$-structure can be distinguished by looking at its restriction to a generator of $H^1(M,\Z/2)$ \cite[Proposition 1 Chapter IV]{kirbybook}.
    It follows that when the loops $\alpha,\alpha'\subset M$ are isotopic, as will be the case in Theorem $\ref{Theorem A}$, the Pin$^+$-structure induced by $(M,\varphi)$ on $\partial\nu(\alpha)$ as the boundary of a codimension zero submanifolds \cite[Construction 1.5]{kirby1990pin} coincides with the one on $\partial\nu(\alpha')$. Let $$f=(f^1,f^2)\colon S^2\simtimes S^1\sqcup S^2\simtimes S^1\to S^2\simtimes S^1\sqcup S^2\simtimes S^1,$$ we denote by $\tau\circ f=(\tau\circ f^1,\tau\circ f^2)$ (see Theorem \ref{isotopyclass}). In this setting we have a natural diffeomorphism $$M(\alpha,\alpha';f)\to M(\alpha,\alpha';\tau\circ f),$$ notice that by the construction described in Lemma \ref{costruzione} the two manifolds are each one Pin$^+$-cobordant to one of the possible classes $(M,+\varphi_M)$, $(M,-\varphi_M)$ but not to each other. 
\end{remark} 
\begin{remark}
    We want to outline possible further use of this construction to construct new stably exotic four-manifolds with small Euler characteristic (up to some group theoretic questions related to presentations of a group). Notice that when we construct $M(\alpha,\alpha',f)$ by choosing disjoint locally flat orientation-reversing simple loops $\alpha,\alpha'\subset M$ such that $[\alpha]\neq[\alpha']\in \pi_1(M)$ that represent two distinct generators, they are related through the gluing described in Lemma \ref{costruzione} to $s$, homotopic to the $S^1$ factor in $S^2\simtimes S^1\times [0,1]$. The new presentation has the same number of generators (we have added the $t$ factor and ``lost" one of the original generators) and the new mapping-torus like relation $sts^{-1}t^{-1}$. In many cases, this procedure reduces the deficiency of the given group presentation by one, but the construction of Lemma \ref{costruzione} preserves the Euler characteristic and gives a cobordism, the fundamental tool used by Kasprowski-Powell \cite[Theorem C, Remark 1.3]{kasprowskipowell} (following Kreck \cite{kreck2006some}) to distinguish smooth structures. It follows that it is interesting to apply the procedure of Lemma \ref{costruzione} to the manifolds in their family with two (or more) generators that are not related. Since Lemma \ref{costruzione} also yields a Pin$^+$-cobordism when the starting manifolds admit a Pin$^+$-structure a similar procedure could be applied to the manifolds described by Torres \cite{torres2025exotic}. 
\end{remark}

One last interesting feature of the gluing described in (\ref{svk}) is that it can be reverted, in the case where the starting loops $\alpha,\alpha'\subset M$ chosen in \textbf{S.1} are homotopic. Let $A=S^1\times [0,1]$ be the trace of the homotopy between $\alpha,\alpha'$, notice that $A$ can be assumed to be an embedded annulus \cite[Example 4.1.3]{gompf19994}. In the gluing $$M(\alpha,\alpha';f)= M_i^0\cup_{f} ((S^2\simtimes S^1)\times[0,1]),$$ the boundary of $A$ is identified with two loops in the boundary of $(S^2\simtimes S^1)\times[0,1]$, and a Klein bottle is described, $$A\cup (\{*\}\times S^1\times[0,1])=Kb\subset M(\alpha,\alpha';f).$$ The explicit nature of the gluing produces a two-sphere $$\Sigma\mathrel{\mathop:}=S^2\times\{q\}\times\{t\}, q\in S^1 \text{ and } t\in[0,1],$$ that intersects transversely in one point the $Kb$. In our context this two-sphere will serve a dual task. Its presence implies that the meridian of $Kb$ is null-homotopic in $M(\alpha,\alpha';f)\setminus\nu(Kb)$, hence $$\pi_1(M(\alpha,\alpha';f))\cong \pi_1(M(\alpha,\alpha';f)\setminus\nu(Kb)),$$ and the construction could lead to some interesting feedstock manifolds to study under operations similar to the ones explored in \cite{cs_cobordism} and \cite[Section 2.1]{akbulut_fs}. Secondly the sphere $\Sigma$, and the bundle structure of $Kb\cong S^1\simtimes S^1$, allow us to pin down the submanifolds needed to get back to our original $M$,  $$\Sigma_{\text{belt}}\mathrel{\mathop:}=(S^2\simtimes S^1\times\{t_0\}),$$ that intersects $Kb$ in a circle. Cutting along $\Sigma_{\text{belt}}$ leaves us with $$M^0=M\setminus(\nu(\alpha)\sqcup\nu(\alpha')),$$ filling up the boundary components with two copies of $D^3\simtimes S^1$ holds the original $M$ \cite[Theorem 3.1]{miller}. 
\subsection{Circle sums}\label{sec:circlesum} Since it is a key tool in the quest to obtain the desired manifolds of Theorem \ref{Theorem A2} and \ref{Theorem B}, we recall the construction introduced in \cite{bais} and \cite{torres2026topologicalclassificationcertainnonorientable}. We consider the total space of the non-orientable three-sphere
bundle over the circle $S^3\simtimes S^1$, and let $\gamma\subset S^3\simtimes S^1$ be a locally flat simple closed loop that generates the fundamental group $\pi_1(S^3\simtimes S^1)\cong\Z$. Let $p>1$ be an odd integer and denote by $\gamma^p\subset S^3\simtimes S^1$ the locally flat orientation-reversing simple closed loop
that represents the homology class $[\gamma]^p\in H_1(S^3\simtimes S^1;\Z)$. Construct the compact topological non-orientable four-manifold \begin{equation}\label{np}
 N_p=(S^3\simtimes S^1)\setminus\nu(\gamma^p),
\end{equation}
with boundary $\partial N_p=S^2\simtimes S^1$.

\begin{theorem}[Bais-Torres {\cite[Theorem B]{bais}}]\label{homeoprocedure}
    Let $M_1$ and $M_2$ be topological four-manifolds and let $$\{\alpha_i\subset M_i:i=1,2\}$$ be locally flat orientation-reversing simple loops. Suppose there is a homeomorphism $$M_1\setminus\nu(\alpha_{1})\to M_2\setminus\nu(\alpha_{2})$$ and define \begin{equation}\label{mapr}M_i(p):=(M_i\setminus\nu(\alpha_i))\cup N_p
    \end{equation}for $i=1,2$, where $N_p$ is the compact four-manifold defined in \ref{np}. The four-manifolds $M_1(p)$ and $M_2(p)$ are homeomorphic. 
\end{theorem}

In the following theorem we summarize the other key features of the procedure that will be needed in the assemblage of the manifolds treated in Theorem \ref{Theorem A2} and Theorem \ref{Theorem B}. 

\begin{theorem}\label{riassunto}\

$\bullet$  \textup{(Bais-Torres) \cite[Lemma 31]{bais}}. The closed non-orientable four-manifold $M(p)$ described in \ref{mapr} has the following basic topological invariants for every $p\in\Z_{>0}$ odd. 

\begin{enumerate}
    \item The Euler characteristic is $\chi(M(p))=\chi(M)$.
    \item The second Stiefel-Whitney class $w_2(M(p))=0$ if and only if $w_2(M)=0$.
\end{enumerate}

$\bullet$  \textup{(Bais-Torres) \cite[Proposition 35]{bais}}. Let $M$ be a closed four-manifold with first cohomology group $H^1(M; \Z/2) = \Z/2$ and that admits a $\Pin^+$-structure $(M, \varphi_M)$. Let $\alpha_M\subset M$ be an orientation-reversing simple loop such that $[\alpha_M]\neq 0\in H_1(M; \Z)$ and build\begin{equation}\label{Construct Cob}M(p)\mathrel{\mathop:}= (M\setminus \nu(\alpha_M))\cup N_p.\end{equation}

There are $\Pin^+$-structures\begin{center}$(S^3\simtimes S^1, \varphi)$ and $(M(p), \varphi_{M(p)})$\end{center}along with a $\Pin^+$-cobordism\begin{equation}\label{Pin Cobordism N}((W, \varphi_W); (M, \varphi_M)\sqcup(S^3\simtimes S^1, \varphi), (M(p), \varphi_{M(p)})).\end{equation}In particular, the identity\begin{equation}\label{Identity Cobordism}[(M, \varphi_M)] = [(M(p), \varphi_{M(p)})]\in \Omega^{\Pin^+}_4\end{equation}holds.   
\end{theorem}
The last cut-and-paste procedure that we want to describe has already been extensively used to produce exotic non-orientable manifolds, see \cite{cappell1976some, HKT,torres2017smooth}. Let $M_1, M_2$ be non-orientable four-manifolds that admit a Pin$^+$-structure. Consider orientation-reversing loops $\gamma_i\subset M_i$ for $i=1,2$. Fix an induced Pin$^+$-structure on the tubular neighbourhood $\nu(\gamma_i)$ and consider \begin{equation}
    i_i:D^3\simtimes S^1\to M_i
\end{equation}
smooth embeddings of $\gamma_i$ such that $i_1$ preserves the Pin$^+$-structure and $i_2$ reverses it. The smooth manifold 

\begin{equation}\label{circlesumeq}
    M_1\#_{S^1}M_2\mathrel{\mathop:}=(M_1\setminus i_1(D^3\simtimes S^1))\cup_f(M_2\setminus i_2(D^3\simtimes S^1))
\end{equation}
is the circle sum of $M_1$ and $M_2$.
    Notice that the circle sum operation is not unique and the resulting manifold depends on the choice of the isotopy class of the gluing map $f$ (see Theorem \ref{isotopyclass}). The key observation for our goals is the following. 

    \begin{proposition}[{\cite[Section 7]{stolz1988exotic}},{\cite[Proposition 1]{torres2017smooth}}]\label{circlesumbordism}
        Let $(M_i,\varphi_i)$ be a closed non-orientable manifold with \textnormal{Pin}$^+$-structure $\varphi_i$ for $i=1,2$. The circle sum $M_1\#_{S^1}M_2$ is \textnormal{Pin}$^+$-bordant to the disjoint union $(M_1,\varphi_1)\sqcup(M_2,\varphi_2)$.
    \end{proposition}

\begin{remark}
   Notice that the first procedure described in Section \ref{sec:circlesum} is a special case of the circle sum construction where $M_2=S^3\simtimes S^1$.
\end{remark}

\subsection{Feedstock four-manifolds}\label{bb} We will need different building blocks to produce the manifolds that are the content of Theorem \ref{Theorem A} and Theorem \ref{Theorem A2}. We give a brief description of such building blocks and focus on their features that will be of relevance in the following. 

In order to build the manifolds described in Theorem \ref{Theorem A} and Theorem \ref{Theorem A2}, our cornerstone will be the manifolds $A(2p)$ and $B(2p)$ produced by Bais-Torres \cite[Theorem F]{bais}, following the work of Cappell-Shaneson \cite{cappell1976some} on exotic $\R\mathbb{P}^4$s.\\
We summarize here the properties of the manifolds $A(2p)$ and $B(2p)$ constructed in \cite[Section 6]{bais}.

\begin{theorem}[Bais-Torres {\cite[Theorem F]{bais}}]\label{startingblockA} For every odd $p\in\Z_{> 0}$, there are closed, smooth, irreducible four-manifolds $A(2p)$ and $B(2p)$ that satisfy the following properties.
\begin{enumerate}
    \item The Euler characteristics of these manifolds are equal to $$\chi(A(2p))=1=\chi(B(2p))$$ and their fundamental groups are isomorphic to $$\pi_1(A(2p))=\Z/2p=\pi_1(B(2p)).$$
    \item There is a homeomorphism $h\colon A(2p)\to B(2p).$
    \item The orientation 2-cover $\widehat{A}(2p)$ is diffeomorphic to the orientation 2-cover $\widehat{B}(2p).$
    \item The universal covers $\widetilde{A}(2p)$ and $\widetilde{B}(2p)$ are diffeomorphic to the connected sum $(p-1)(S^2\times S^2)$.
    \item There is no diffeomorphism $$A(2p)\#(k-1)(S^2\times S^2) \to B(2p)\#(k-1)(S^2\times S^2)$$ for any $k\in\Z_{> 0}.$
    \item The connected sums $A(2p)\#\CP^2$ and $B(2p)\#\CP^2$ are diffeomorphic.   
\end{enumerate}
    Moreover, there is a non-smoothable closed topological four-manifold $*A(2p)$ that is homotopy equivalent to $A(2p)$.
\end{theorem}

\section{Assemblage and certification of the four-manifolds of theorem \ref{Theorem A} and theorem \ref{Theorem A2}}\label{produccion} 

\subsection{Construction of the four-manifolds of Theorem \ref{Theorem A}}\label{constructA} 
To build our examples, we plug in the four-manifolds of Theorem \ref{startingblockA}, $A(2p)$ and $B(2p)$, homeomorphic through $h\colon A(2p)\to B(2p)$, into the mechanism described in Section \ref{lego}. We choose disjoint simple loops $\{\alpha,\alpha'\}\subset A(2p)$ and $\{\beta,\beta'\}:=\{h(\alpha),h(\alpha')\}\subset B(2p)$ such that $\alpha,\alpha'$ are homotopic to the generator of $\pi_1(A(2p))$. We obtain 
\begin{equation}\label{newmani1}
    (X_{2p},\varphi_{X_{2p}}):=A(2p)(\alpha,\alpha',f_a)
\end{equation} and 
\begin{equation}\label{newmani2}
    (Y_{2p},\varphi_{Y_{2p}}):=B(2p)(\beta,\beta',f_b)
\end{equation} where $f_a$ and $f_b=f_a\circ h^{-1}$ have been chosen so that they extend the Pin$^+$-structure and produce the desired fundamental group presentation \begin{align}\label{fundamentalA}
\pi_1(X_{2p})&=\langle\alpha,\alpha',s,t\;|\;\alpha\alpha'^{-1}, \alpha^{2p},\alpha s,\alpha's,s ts^{-1}t^{-1}\rangle \\
\nonumber &=\langle\alpha,t\;|\;\alpha^{2p},\alpha t\alpha^{-1}t^{-1}\rangle\cong\Z\times\Z/2p
\end{align} as described in Lemma \ref{costruzione}. Notice that for Proposition \ref{homeocondition} and Remark \ref{pinremark} to carry on here with the appropriate choice of Pin$^+$-structures is key that $\alpha,\alpha'\subset A(2p)$ are homotopic loops in a 4-manifold, hence isotopic \cite[Example 4.1.3]{gompf19994}, and a basis for $H^1(A(2p);\Z/2)\cong\Z/2$ \cite[Proposition 1 Chapter IV]{kirbybook}.

Observe that the possible Pin$^+$-structures on $A(2p)$ and $B(2p)$ are extensively discussed in \cite[Proposition 35, Proof of Theorem F]{bais}, the Pin$^+$-cobordism of Lemma \ref{costruzione}, Theorem \ref{Theorem NonOrientable 3} and the computations on the $\eta$-invariant in \cite{bais} imply the following 

\begin{equation}\label{eta1}
    \eta(X_{2p},\pm\varphi_{X_{2p}})=\eta(A(2p),\pm \varphi_{A(2p)})=\pm\frac{1}{8} \;\text{mod}\; 2\Z
\end{equation}

\begin{equation}\label{eta2}
    \eta(Y_{2p},\pm\varphi_{Y_{2p}})=\eta(B(2p),\pm \varphi_{B(2p)})=\pm\frac{7}{8} \;\text{mod}\; 2\Z.
\end{equation}
\begin{remark}
    Since $H^1(X_{2_p};\Z/2)\cong H^1(Y_{2p};\Z/2)\cong \Z/2\times \Z/2$, both manifolds admit four distinct Pin$^+$-structures. These \textsl{additional} Pin$^+$-structures can be distinguished by looking at their restriction to the loop $t\subset M(\alpha,\alpha';f)$ coming from the $I$- factor of $S^2\simtimes S^1\times I$ \cite[ Proposition 1 Chapter IV]{kirbybook}. Starting from $M(\alpha,\alpha';f)$ and running the cobordism argument of Lemma \ref{costruzione} backwards identifies the remaining Pin$^+$
-structures of $W$ and $W'$ from Lemma, via the natural inclusions $X_{2p}\hookrightarrow W$ and $Y_{2p}\hookrightarrow W'$ which induce isomorphism of $H^1(-;\Z/2)$. Notice that the choice of the restriction of the Pin$^+$-structure to the loop $t$ is immaterial when we attach the boundary of the cylinder $S^2\simtimes S^1\times I$ to the loops $\alpha,\alpha'$.  In our case of Theorem \ref{Theorem A}, while the $\eta$-invariant does not distinguish individual structures within each manifold, the sets of $\eta$-values for $X_{2p}$ and $Y_{2p}$ are completely disjoint modulo $2\Z$, \begin{align}\label{pinextension}
        &\eta(X_{2p},\pm\varphi_{X_{2p}}^i)=\pm\frac{1}{8} \;\text{mod}\; 2\Z\;\text{and}\;\\&\eta(Y_{2p},\pm\varphi_{Y_{2p}}^i)=\pm\frac{7}{8}\;\text{mod}\; 2\Z\;\textnormal{for }i\in\{1,2,3,4\},
    \end{align}which is sufficient to apply Theorem \ref{torresexotic}. 
\end{remark}
Now, we will consider the orientation 2-covers $\widehat{A}(2p)$ and $\widehat{B}(2p)$. Denote the preimages in the coverings of the loops $\alpha,\alpha',\beta,\beta'$ with $\widehat{\alpha},\widehat{\alpha}',\widehat{\beta},\widehat{\beta}'$. Since these lifted loops are submanifolds of an oriented manifold their tubular neighbourhoods are all diffeomorphic to $D^3\times S^1$. By considering the lift of $f\colon S^2\simtimes S^1\sqcup S^2\simtimes S^1\to S^2\simtimes S^1\sqcup S^2\simtimes S^1$ to a $\widehat{f}\colon S^2\times S^1\sqcup S^2\times S^1\to S^2\times S^1\sqcup S^2\times S^1$ \cite[Proof of Theorem 2.3]{kim1990diffeotopy} we obtain the degree two coverings of $X_{2p}$ and $Y_{2p}$. 

\begin{equation}\label{orientablecovering1}
    \widehat{X}_{2p}= (\widehat{A}(2p)\setminus(\nu(\widehat{\alpha})\sqcup\nu(\widehat{\alpha}')))\cup_{\widehat{f}}(S^2\times S^1\times I).
\end{equation}

\begin{equation}\label{orientablecovering2}
    \widehat{Y}_{2p}= (\widehat{B}(2p)\setminus(\nu(\widehat{\beta})\sqcup\nu(\widehat{\beta}')))\cup_{\widehat{f}}(S^2\times S^1\times I).
\end{equation}

\subsection{Proof of Theorem \ref{Theorem A}}\label{certificateA}
The two four-manifolds of Theorem \ref{Theorem A} were constructed in (\ref{newmani1}), (\ref{newmani2}). The homeomorphism class of these four-manifolds coincides, as the natural homeomorphism is the one obtained by extending  \begin{equation}
    h:A(2p)\setminus(\nu(\alpha)\sqcup\nu(\alpha'))\rightarrow B(2p)\setminus(\nu(\beta)\sqcup\nu(\beta')),
\end{equation} the canonical restriction of the homeomorphism from Theorem \ref{startingblockA} as described in Proposition \ref{homeocondition}. The fundamental group and Euler characteristic of $X_{2p}$ are as desired, as a consequence of Lemma \ref{costruzione}, Theorem \ref{startingblockA} and the computation of (\ref{fundamentalA}). Item (3) is a direct consequence of Theorem \ref{torresexotic} and the computations of (\ref{pinextension}). Item (4) follows from the constructions \ref{orientablecovering1}, \ref{orientablecovering2} and the diffeomorphism of Item (3) from \ref{startingblockA}. The diffeomorphism of Item (5) is inherited by construction from the diffeomorphism described in Theorem \ref{startingblockA} and by performing the connected sum choosing disks disjoint from the loops $\{\alpha,\alpha',\beta,\beta'\}$; \begin{align*}
    X_{2p}\#\CP^2&:= A(2p)(\alpha,\alpha';f_a)\#\CP^2\cong (A(2p)\#\CP^2)(\alpha,\alpha';f_a)\\
    &\cong(B(2p)\#\CP^2)(\beta,\beta';f_b)\cong B(2p)(\beta,\beta';f_b)\#\CP^2\\&=: Y_{2p}\#\CP^2.
\end{align*} \qed 

\begin{remark}\label{rmk1}
    Our manifolds have Euler characteristic equal to 1 and fundamental group $\pi_1(X_{2p})=\Z\times\Z/2p$. The presentation $$\Z\times\Z/2p=\langle x,y\;|\;[x,y],y^{2p}\rangle $$ realizes our group and has a deficiency of 0 (cf. Appendix \ref{deficiency}).
\end{remark}

Except for $p=1$, the above-mentioned manifolds minimize the Euler characteristic in their fundamental group class. 

\begin{proposition}\label{minimize}
    Let $p\in\Z_{>1}$ be an odd integer. Let $M$ be a closed, non-orientable four-manifold such that $\pi_1(M)=\Z\times\Z/2p$, then $\chi(M)\geq 1$.
\end{proposition}
\begin{proof}
    Let $G=\Z\times\Z/2p$. Consider the orientation double cover $$\widehat{M}\overset{2:1}\to M,$$ let $K=\pi_1(\widehat{M})$, this is an index-two subgroup of $\pi_1(M)$, hence it is isomorphic to either $\Z\times\Z/2p$ or $\Z\times\Z/p$ and in either case we have a natural inclusion $i\colon K\to G$. Let $q$ be a prime number such that $q\mid p$. Now consider the following commutative (up to homotopy) diagram
        \[
\begin{tikzcd}[row sep=large, column sep=large]
  \widehat M \arrow[r, "\pi"] \arrow[d, "\widehat c"'] 
  & M \arrow[d, "c"] \\
  BK \arrow[r, "Bi"'] 
  & BG
\end{tikzcd}
\] where classifying maps are induced by the respective fundamental-group isomorphisms and let $\alpha=\widehat{c}_*[\widehat{M]}\in H_4(K;\Z)$. 
Observe that the non-trivial deck transformation $$\tau\colon\widehat{M}\to \widehat{M},$$ is an orientation reversing diffeomorphism of $\widehat{M}$ and let it be represented by an element $g\in G\setminus i(K)$. We want to show that $\widehat{c}\circ \tau\simeq \widehat{c}$. Up to change of basepoint the action that $\tau$ induces on $\pi_1(\widehat{M})\cong K$ is conjugation by $g$: $$\tau_*(k)=gkg^{-1}, \textnormal{ for every } k\in K.$$ Given that $G$ is abelian this is the identity on $K$, hence $$\widehat{c}\circ \tau\simeq \widehat{c}.$$ Now we have $$\alpha=\widehat{c}_*[\widehat{M}]=(\widehat{c}\circ\tau)_*[\widehat{M}]=\widehat{c}_*\circ\tau_*[\widehat{M}]=-\widehat{c}_*[\widehat{M}]=-\alpha,$$  
where the penultimate equality follows from the orientation reversing nature of the deck transformation $\tau$. As a consequence, we have $2\alpha=0$, and given that $q$ is odd $$q\alpha=\alpha+(q-1)\alpha=\alpha.$$ We have proven that $\alpha$ is a multiple of $q$, it follows from \cite[Theorem 4.6]{Kirk} (cf. \cite[Lemma 1.5]{kreckluck}) \begin{align}\label{cohomologyeq}
    \chi(\widehat{M})&\geq 2-2\textnormal{dim}_{\Z/q}H^1(K;\Z/q)+2\textnormal{dim}_{\Z/q}H^2(K;\Z/q)\\
    &=2-2\cdot2+2\cdot2=2\nonumber.
\end{align}
See Appendix \ref{deficiency} for the computation of the cohomology groups appearing in \eqref{cohomologyeq}. It is a straightforward consequence of the multiplicativity of Euler characteristic on finite degree coverings that $$\chi(M)\geq 1.$$

\end{proof}

\subsection{Construction of the four-manifolds of Theorem \ref{Theorem A2}}\label{sec:a2}
We iterate the construction devised by Bais-Torres in \cite[Section 6]{bais} and described in Section \ref{sec:circlesum}, starting from the manifold $A(2p)$, $B(2p)$ of Theorem \ref{startingblockA}.
Let $p\in\Z_{>1}$ odd and $\alpha\subset A(2p)$, be a simple loop whose homotopy
class generates $\pi_1(A(2p))=\Z/2p$. Denote
by $\alpha_p\subset A(2p)$ the simple loop that represents the class $[\alpha]^p\in \pi_1(A(2p))$, and by $\beta_p=h(\alpha_p)$ the image of $\alpha_p$ under the homeomorphism $h\colon A(2p)\to B(2p)$ of Theorem \ref{startingblockA}. Note that this element has order $2$ and its orientation character is $w(\alpha_p)=1$. Fix an odd $q_1\in\Z_{>0}$ and consider the compact four-manifold $N_{q_1}$ of \eqref{np}. Use these data as input in the Bais-Torres procedure to assemble the non-orientable four-manifolds \begin{align}\label{noq1}
    A(2p)(q_1):= (A(2p)\setminus\nu(\alpha_p))\cup N_{q_1}
\end{align}
\begin{align}\label{nobq1}
    B(2p)(q_1):= (B(2p)\setminus\nu(\beta_p))\cup N_{q_1}
\end{align}

The Seifert-van Kampen theorem holds \begin{align*}
    \pi_1(A(2p)(q_1))&=\langle\alpha,\gamma \;|\;\alpha^{2p},\alpha^p\gamma^{-q_1}\rangle\\\notag &=\langle\alpha,\gamma,\zeta,\mu \;|\;\alpha^{2p},\alpha^p\gamma^{-q_1}, \alpha^p\zeta^{-1},\alpha^2\mu^{-1}\rangle\\\notag&=\langle\alpha,\gamma,\zeta,\mu \;|\;\alpha^{2p},\alpha^p\gamma^{-q_1},\alpha^p\zeta^{-1}, \alpha^2\mu^{-1},\zeta^2, [\zeta,\alpha],[\zeta,\gamma]\rangle\\
    \notag &=\langle\alpha,\gamma,\zeta,\mu,\nu \;|\;\alpha^{2p},\alpha^p\gamma^{-q_1},\alpha^p\zeta^{-1}, \alpha^2\mu^{-1},\zeta^2, [\zeta,\alpha],[\zeta,\gamma], \gamma^2\nu^{-1}\rangle \\
    \notag &=\langle\zeta,\mu,\nu|\mu^p,\nu^q_1,\zeta^2,[\zeta,\mu],[\zeta,\nu]\rangle=(\Z/2)\times(\Z/p*\Z/q_1).
\end{align*} The last simplifications in the presentation are motivated by the following relations, \begin{align*}
    \alpha&=\alpha^{2p+1}=\alpha^{p}\alpha^{p+1}=\zeta\mu^{\frac{p+1}{2}},\\
    \gamma&=\gamma^{2q_1+1}=\gamma^{q_1}\gamma^{q_1+1}=\zeta\nu^{\frac{q_1+1}{2}}.
\end{align*}Let $\zeta_a\subset A(2p)(q_1)$ be the simple closed loop that represents the class of $[\zeta]\in\pi_1(A(2p)(q_1))$. Notice that $\zeta_a$ is a non-orientable loop of order two just as $\alpha_p\subset A(2p)$. Fix an odd $q_2\in\Z_{>0}$, consider the compact four-manifold  $N_{q_2}$ of \eqref{np} and plug in these new data in the Bais-Torres procedure to obtain 
\begin{align}
    A(2p)(q_1,q_2)=(A(2p)(q_1)\setminus\nu(\zeta_a))\cup N_{q_2}.
\end{align}
    The Seifert-van Kampen theorem holds
    \begin{align*}
        \pi_1(A(2p)(q_1,q_2))&=\langle\zeta,\mu,\nu,\gamma\;|\;\mu^p,\nu^{q_1},\zeta^2,[\zeta,\mu],[\zeta,\nu],\zeta\gamma^{-q_2}\rangle \\
        &=\langle\zeta,\mu,\nu,\gamma,\lambda\;|\;\mu^p,\nu^{q_1},\zeta^2,[\zeta,\mu],[\zeta,\nu],\zeta\gamma^{-q_2},\gamma^2\lambda^{-1}\rangle \\
        &=\langle\zeta,\mu,\nu,\lambda\;|\;\mu^p,\nu^{q_1},\lambda^{q_2},\zeta^2, [\zeta,\mu],[\zeta,\nu],[\zeta,\lambda]\rangle\\&=((\Z/2)\times(\Z/p*\Z/q_1*\Z/q_2)).
    \end{align*}
    Here we are using a similar algebraic equality as before \[\gamma=\gamma^{2q_2+1}=\gamma^{q_2}\gamma^{q_2+1}=\zeta\lambda^{\frac{q_2+1}{2}}.\]
    Notice that in the presentation of the fundamental group we still have the distinct non-orientable element of order two $\zeta$ and the new simplifications did not involve any other old generators. This suggests iterating the construction as follows. Suppose we have constructed $A(2p)(q_1,\dots,q_{m-1})$, fix an odd $q_m\in\Z_{>0}$ consider the compact four-manifold $N_{q_m}$ and run once again the Bais-Torres procedure to obtain $A(2p)(q_1,\dots,q_m)$, a closed non-orientable four-manifold with fundamental group \[\pi_1(A(2p)(q_1,\dots,q_m))=((\Z/2)\times(\Z/p*\Z/q_1*\Z/q_2*\dots*\Z/q_m)).\] 

    \subsection{Proof of Theorem \ref{Theorem A2}}\label{sec:proofA2}

   We have assembled the manifolds $A(2p)(q_1)$ and $B(2p)(q_1)$ in \eqref{noq1} and \eqref{nobq1}, the homeomorphism $$h\colon A(2p)\to B(2p)$$ extends to a homeomorphism $$h_{q_1}\colon A(2p)(q_1)\to B(2p)(q_1)$$ by Theorem \ref{homeoprocedure}. This  allows us to run the iterative procedure described in Section \ref{sec:a2} and obtain homeomorphism between the constructed manifolds at each step. The computations of the fundamental groups have been described during the construction of the manifolds in Section \ref{sec:a2}, the iterative nature of the construction, and Theorem \ref{riassunto}, showing
   \begin{align*}
       \chi(A(2p))&=\chi(A(2p)(q_1))=\dots=\chi(A(2p)(q_1,\dots,q_m))=1, \\
       w_2(A(2p))&=w_2(A(2p)(q_1))=\dots=\   w_2(A(2p)(q_1,\dots,q_m))=0.
   \end{align*}

   The extension of Pin$^+$-structures and the cobordism described in Theorem \ref{riassunto} permit to conclude the stable exoticity of the third Item of Theorem \ref{Theorem A2}. For every $q\in \Z_{>1}$ consider the orientation 2-cover \[S_q\to N_q\] described by Bais-Torres in \cite[Lemma 7]{bais} and the diffeomorphism between the orientation 2-covers $\widehat{A}(2p)$ and $\widehat{B}(2p)$ of Theorem \ref{startingblockA}. At the first step of our construction, the loops $\alpha_p \in A(2p)$ and $\beta_p\subset B(2p)$ have as preimage in the 2-covers null-homotopic orientable loops which are glued to $S_q$. From \cite[Lemma 7]{bais} we get the following diffeomorphism, 
   \[\widehat{A(2p)(q_1)}\cong\widehat{A(2p)}\#\Sigma_{q_1}\cong\widehat{B(2p)}\#\Sigma_{q_1}\cong\widehat{B(2p)(q_1)},\]
 where $\Sigma_{q_1}$ is a rational homology sphere with $\pi_1(\Sigma_{q_1})=\Z/q_1$. Playing the same game with loops $\zeta_a$ and $\zeta_b$ allows us to iterate the diffeomorphism of the orientation 2-covers. For the last item, we perform the connected sum along disks disjoint from the surgery loops and the result follows again from the diffeomorphism of Theorem \ref{startingblockA}. \qed
 
\begin{proposition}\label{secondminimality}
    Let $p,q_1,\dots,q_m\in\Z_{>1}$ be odd integers. Let $M$ be a closed, non-orientable four-manifold such that $$\pi_1(M)=\Z/2\times(\Z/p*\Z/q_1*\dots*\Z/q_m),$$ then $\chi(M)\geq 1$.
\end{proposition}
\begin{proof}
    Let $G=\Z/2\times(\Z/p*\Z/q_1*\dots*\Z/q_m)$. Let $M$ be a connected, compact, non-orientable four-manifold with $\pi_1(M)=G$. Poincaré duality with $\Z/2$ coefficients implies 
    $$\chi(M)=2-2b_1(M;\Z/2)+b_2(M;\Z/2).$$ 
    This fact, Hopf's theorem \cite[Theorem 5.2]{Brown} and the computations of equation (\ref{lastkunnet}) carried out in Appendix \ref{deficiency} imply \[\chi(M)\geq2-2b_1(G;\Z/2)+b_2(G;\Z/2)=1,\] hence that the manifolds produced in Theorem \ref{Theorem A2} minimize the Euler characteristic among closed non-orientable four-manifolds in their fundamental group class.
\end{proof}
\section{On the non-orientable Hausmann-Weinberger invariant}\label{HWinvariant}
We continue with some investigation of the non-orientable Hausmann-Weinberger invariant. Recall the definition from Torres \cite[Definition D]{torres2026topologicalclassificationcertainnonorientable}.

\begin{definition}
    Let $G$ be a finitely presented group. The non-orientable
Hausmann-Weinberger invariant $q(G)$ of $G$ is the minimal value of the
Euler characteristic of a closed non-orientable topological four-manifold $M$
with fundamental group $G$.
\end{definition} First, we reap the rewards of the various cut-and-paste mechanism scattered in the note, after that we rely on group cohomology computations to show that certain simple product manifolds are minimizers.
\subsection{Proof of Theorem \ref{Theorem B}}\label{proofofB}
The manifolds of Theorem \ref{Theorem A} and Theorem \ref{Theorem A2} have exactly the fundamental group and the Euler characteristic described in Items (1) and (2). Minimality is the content of Proposition \ref{minimize} and Proposition \ref{secondminimality}. \\
For the third item, start by considering 
\begin{equation}\label{vu}
    V=\mathbb{RP}^2\simtimes \mathbb{RP}^2
\end{equation} 
Notice that $\chi(V)=1$ and fix a presentation of the Klein group $$\pi_1(V)=\langle  \alpha,\beta\;|\;\alpha^2,\beta^2,[\alpha,\beta]\rangle.$$
Consider a locally flat orientation-reversing simple loop $\alpha\subset V$ that is a generator for $\pi_1(V)$ and remove the tubular neighbourhood of the loop, to obtain a compact four-manifolds with boundary that can be plugged in the procedure described in Section \ref{sec:circlesum}. Carrying on with the following steps, we obtain \begin{equation}\label{mp}
    V(p):=V\setminus \nu(\alpha)\cup N_p.
\end{equation}

Applying the Seifert-van Kampen theorem to the decomposition of \eqref{mp} we get, 

\begin{align*}\label{presentazione}
        \pi_1(V(p))&=\langle \alpha,\beta,\gamma \;|\; \alpha^2,\beta^2,[\alpha,\beta],\;\alpha\gamma^{-p}\rangle \\
        &=\langle\alpha,\beta,\gamma,\lambda\;|\;\alpha^2,\beta^2,[\alpha,\beta],\;\alpha\gamma^{-p},\gamma^{2}\lambda^{-1}\rangle\\
        &=\langle\alpha,\beta,\lambda\;|\;\alpha^2,\beta^2,\lambda^p,[\alpha,\beta],[\alpha,\lambda]\rangle\\
        &=\Z/2\times(\Z/2*\Z/p).
    \end{align*}. 
Here we have used the same algebraic manipulation of Section \ref{sec:a2}, 
\[\gamma=\gamma^{2p+1}=\gamma^p\gamma^{p+1}=\alpha\lambda^{\frac{p+1}{2}}.\]

With the same construction iterating the gluing and the algebraic manipulation with $N_{q_1},\dots N_{q_m}$ we obtain $V(p,q_1,\dots q_m)$ with fundamental group \[\pi_1(V(p,q_1,\dots q_m))=\Z/2\times(\Z/2*\Z/p*\Z/q_1\dots*\Z/q_m).\]

Theorem \ref{riassunto} allows us to conclude that $$\chi(V(p,q_1,\dots,q_m))=\dots=\chi(V(p))=\chi(V)=1.$$ 

Let $G=\Z/2\times(\Z/2*\Z/p*\Z/q_1\dots*\Z/q_m)$ and let $M$ be a connected, compact, non-orientable four-manifold with $\pi_1(M)=G$. Poincaré duality with $\Z/2$ coefficients implies 
    $$\chi(M)=2-2b_1(M;\Z/2)+b_2(M;\Z/2).$$ 
    This fact, Hopf's theorem \cite[Theorem 5.2]{Brown} and the computations of equation (\ref{bettinumb}) carried out in Appendix \ref{deficiency} imply \[\chi(M)\geq2-2b_1(G;\Z/2)+b_2(G;\Z/2)=1,\] hence that the manifolds $V(p,q_1,\dots q_m)$ produced minimize the Euler characteristic among closed non-orientable four-manifolds in their fundamental group class.

    \begin{remark}
        The group presentation \[G=\langle \alpha,\beta,\gamma \;|\; \alpha^2,\beta^2,[\alpha,\beta],\;\alpha\gamma^{-p}\rangle=\Z/2\times(\Z/2*\Z/p),\] for $p\in\Z_{>1}$ odd, already appeared in the work of the author and Torres \cite{camponovotorres}. There it was only remarked that $G$ is infinite with exponential growth; we did not notice at the time the algebraic manipulation that prove the isomorphism with $\Z/2\times(\Z/2*\Z/p)$.
    \end{remark}

\subsection{Proof of Theorem \ref{Theorem C}}
    We begin by proving Proposition \ref{propcoveringq(G)}.
         Let $C_0<G$ be an infinite cyclic subgroup with finite index. Let \[w:\pi_1(M)\to\Z/2,\] be the orientation character and $H=\ker w$. $M$ is non-orientable, so $[\pi_1(M):H]=2$. Consider the subgroup \[C=C_0\cap H.\] $C$ has finite index in $\pi_1(M)$ and is an infinite subgroup of $C_0$, hence $C\cong\Z$. Consider the finite degree covering \[p\colon\widehat{M}\to M\] associated to $C$. Since $C\subset H$, the cover space is an orientable manifold and $\pi_1(\widehat{M})\cong\Z$, so \[b_1(\widehat{M})=b_3(\widehat{M})=1.\]
    Given that $\widehat{M}$ is oriented and connected, 
    \[b_0(\widehat{M})=b_4(\widehat{M})=1.\]
    Therefore, 
    \[\chi(\widehat{M})=1-1+b_2(\widehat{M})\geq0,\]
    from the multiplicativity of the Euler characteristic under finite covering, \[\chi(\widehat{M})=[\pi_1(M):C]\chi(M),\] hence \[\chi(M)\geq0.\]     
    The proof of Theorem \ref{Theorem C} is now just an application of Proposition \ref{propcoveringq(G)} and the observation \[\chi(M_f)=\chi(M)\chi(S^1)=0.\] \qed
   
    \subsection{Proof of Theorem \ref{Theorem D}}
    
   Let $M$ be a closed connected, compact, non-orientable four-manifold. Poincaré duality with $\Z/2$ coefficients implies 
    $$\chi(M)=2-2b_1(M;\Z/2)+b_2(M;\Z/2).$$ Let $G=\pi_1(M)$, then by Hopf's theorem \cite[Theorem 5.2]{Brown} \[\chi(M)\geq2-2b_1(G;\Z/2)+b_2(G;\Z/2).\] 
Therefore, via the computation in \eqref{bettisurface} 
    \[q(\Gamma_g)\geq2-2(2g+1)+(2g+2)=2-2g.\]
    Then it follows that the four-manifold  $\Sigma_g\times\RP^2$ realizes $q(\Gamma_g)$. \\
    For the second part of the statement, consider $L(p,q)\simtimes S^1$ assembled as described in Remark \ref{nonorientablelens}. Now choose a loop $\alpha\subset L(p,q)\simtimes S^1$ representing the section of the mapping torus and a loop $\gamma\subset\RP^4$ generating $\pi_1(\RP^4)$. Let $r\in\Z_{>1}$ be an odd integer and perform the circle sum operation described in \ref{sec:circlesum} with simple closed loops representing $\alpha^r$ and $\gamma$ to obtain \[M_{p,q,r}=(L(p,q)\simtimes S^1)\#_{S^1}\RP^4.\]
    The Seifert-van Kampen theorem holds 
\begin{align*}
    \pi_1(M_{p,q,r})&=\langle \alpha,\beta,\gamma\;|\;\beta^p,\alpha\beta\alpha^{-1}\beta^{-q},\gamma^2,\alpha^r\gamma^{-1}\rangle \\ &=\langle \alpha,\beta\;|\;\beta^p,\alpha\beta\alpha^{-1}\beta^{-q},\alpha^{2r}\rangle \\
    &=\langle \alpha,\beta\;|\;\beta^2,\alpha^{2r},[\alpha,\beta]\rangle=\Z/2\times\Z/2r.
\end{align*}
The last step follows from the $q^2\equiv -1 \pmod{p}$ and $r$ being odd: \[\beta=\alpha^{2r}\beta\alpha^{-2r}=\beta^{q^{2r}}=(\beta^{-1})^r=\beta^{-1}.\]
For the Euler characteristic we have \[\chi(M_{p,q,r})=\chi(L(p,q)\simtimes S^1)+\chi(\RP^4)=0+1=1\] and once again by Hopf's Theorem \cite[Theorem 5.2]{Brown} \[\chi(M_{p,q,r})\geq2-2b_1(\Gamma_r;\Z/2)+b_2(\Gamma_r;\Z/2),\]
from the computation in \eqref{ultimissimissima}\[q(\Gamma_r)\geq 2-2(2)+3=1.\] \qed
    \appendix
    \section{On group deficiency and some (co)homological calculations.}\label{deficiency}
For a finitely presented group $G$, let $\textnormal{def}(G)$ be the deficiency of $G$, that is the maximum over
all presentations of the number of generators minus the number of relators. Here we present some calculations with group (co)homology that are paramount for the proofs in this note. See \cite{Brown} for further details on the tools used. One of the key results that we will make use of is the following. 

\begin{lemma}[{\cite[Lemma 3.2]{Jin}}]\label{rose}
    For any finitely presented group $G$ and any field $F$,

    \begin{equation}
        \textnormal{def}(G)\leq b_1(G;F)-b_2(G;F).
    \end{equation}
\end{lemma}

Recall that $S^1=K(\Z,1)$ (cf. \cite[Proposition 4.1]{Brown}), using the universal coefficient theorem we have 

\begin{equation}   
H_i(\Z;\Z/2)
\cong
\begin{cases}
\Z/2, & i=0,1,\\
0, & \text{otherwise}.
\end{cases}
\end{equation}

By the definition given in \cite{Brown}, we get that for every $p\in\Z_{>0}$, 

\begin{equation}\label{mod2}
    H_i(\Z/2p;\Z/2)\cong \Z/2,\ \text{for every }i>0.
\end{equation}

By the K\"unneth formula \cite[Chapter 3]{HatcherAT}, we obtain 
\begin{align*}
    &H_1(\Z\times\Z/2p;\Z/2)\cong H_1(\Z;\Z/2)\oplus H_1(\Z/2p;\Z/2), \\
   &H_2(\Z\times\Z/2p;\Z/2)\cong H_1(\Z;\Z/2)\otimes H_1(\Z/2p;\Z/2)\oplus H_2(\Z/2p;\Z/2).
\end{align*}
From this we get $$b_1(\Z\times\Z/2p;\Z/2)=2,$$
and $$b_2(\Z\times\Z/2p;\Z/2)=1+1=2.$$
Therefore, via Lemma \ref{rose} we get the following bound 
$$\text{def}(\Z\times\Z/2p)\leq b_1(\Z\times\Z/2p;\Z/2)-b_2(\Z\times\Z/2p;\Z/2)=0.$$
The presentations of Remark \ref{rmk1}, allows us to conclude that $$\text{def}(\Z\times\Z/2p)=0, \text{ for }p\in\Z_{>0}.$$

   We now carry out the computations needed for the proof of Proposition \ref{minimize}. Let $p>1$ be an odd integer, let $m\in\{p,2p\}$ and $q$ a prime such that $q\mid p$. Reasoning as above we have: \begin{equation}   
H^i(\Z;\Z/q)
\cong
\begin{cases}
\Z/q, & i=0,1,\\
0, & \text{otherwise}.
\end{cases}
\end{equation}

For the finite cyclic group, given that \(q\mid m\), we have the following, from the standard periodic resolution \cite[Chapter 3]{Brown}.

$$H^i(\Z/m;\Z/q)\cong \Z/q, \text{ for every }i\geq0.$$ From the K\"unneth formula
$$H^1(\Z\times\Z/m;\Z/q)\cong (\Z/q)^2,$$
and 
$$H^2(\Z\times\Z/m;\Z/q)\cong (\Z/q)^2,$$ as claimed in Proposition \ref{minimize}.

To compute Betti numbers of the groups described in Theorem \ref{Theorem B}  we will need the following. \begin{lemma}[{\cite[Corollary 7.7]{Brown}}]\label{MV}
    Let $G_1,G_2$ and $A$ be groups with monomorphisms $A\hookrightarrow G_1,$ $A\hookrightarrow G_2$, we denote by $G=G_1*_AG_2$ the amalgamated product. There is a Mayer-Vietoris sequence \begin{equation*}
        \dots\to H_n(A)\to H_n(G_1)\oplus H_n(G_2)\to H_n(G)\to H_{n-1}(A)\to\dots
    \end{equation*}
        
\end{lemma} 

Now we let $H=\Z/p*\Z/q_1*\dots*\Z/q_m$ for $p,q_1,\dots q_m\in \Z_{>1}$ odd integers and $G=\Z/2\times H$, from here we have \cite[Example 1B.5]{HatcherAT} \[K(G,1)=K(\Z/2,1)\times K(H,1).\]\ For $p$ odd we have \cite[Chapter 3]{Brown} \begin{equation}   
H_i(\Z/p;\Z/2)
\cong
\begin{cases}
\Z/2, & i=0,\\
0, & \text{otherwise}.
\end{cases}
\end{equation}

From Lemma \ref{MV} we have \[H_n(G_1*G_2;\Z/2)\cong H_n(G_1;\Z/2)\oplus H_n(G_2;\Z/2),\;n\geq1.\]

By the K\"unneth formula, and recalling $K(\Z/2,1)=\RP^\infty$,

\begin{equation}\label{lastkunnet}
    H_i(G;\Z/2)\cong \Z/2, i\in\{1,2\},
\end{equation}
as claimed in Proposition \ref{secondminimality}. From Lemma \ref{rose} and the group presentation of \ref{sec:a2} we have $\text{def}(G)=0$.

Now we let $H=\Z/2*\Z/q_1*\dots*\Z/q_m$ for $q_1,\dots q_m\in \Z_{>1}$ odd integers and $G=\Z/2\times H$. The proof follows verbatim as before, therefore \begin{equation}\label{bettinumb} H_1(G;\Z/2)\cong (\Z/2)^{2},\;H_2(G;\Z/2)\cong(\Z/2)^3,
\end{equation}
as claimed in the proof of Theorem \ref{Theorem B} (Section \ref{proofofB}). 

Let $\Sigma_g$ be a closed orientable surface of genus $g\geq1$. For the group $\Gamma_g=\pi_1(\Sigma_g\times\RP^2)=\pi_1(\Sigma_g)\times \Z/2$ we can compute (cf. \cite[Proposition 4.1]{Brown}) the Betti numbers using \[K(\Gamma_g,1)=\Sigma_g\times\RP^\infty.\]

    From the K\"unneth formula \begin{equation}\label{bettisurface}
        b_1(\Gamma_g;\Z/2)=2g+1,\\
b_2(\Gamma_g;\Z/2)=2g+2
    \end{equation}

    Let $r\in \Z_{>1}$ odd, from the K\"unneth formula and using \eqref{mod2} we have 

    \begin{equation}\label{ultimissimissima} b_1(\Z/2\times\Z/2r;\Z/2)=2,\\b_2(\Z/2\times\Z/2r;\Z/2)=3.
    \end{equation}

\section*{References}
\printbibliography[heading=none]
\end{document}